# A new operator giving integrals and derivatives operators for any order at the same time

### Raoelina Andriambololona

Theoretical Physics Department, Institut National des Sciences et Techniques Nucléaires (INSTN-Madagascar)

instn@moov.mg, raoelinasp@yahoo.fr, jacquelineraoelina@yahoo.fr, raoelina.andriambololona@gmail.com

*Abstract- Let E be the set of integrable and derivable causal functions of $x$ defined on the real interval I from $a$ to infinity, $a$ being real, such $f(a)$ is equal to zero for $x$ lower than or equal to $a$. We give the expression of one operator that yields the integral operator and derivative operators of the function $f$ at any s-order. For $s$ positive integer real number, we obtain the ordinary s-iterated integral of $f$. For $s$ negative integer real number we obtain the $|s|$-order ordinary derivatives of $f$. Any $s$ positive real or positive real part of $s$ complex number corresponds to s-integral operator of $f$. Any $s$ negative real number or negative real part of $s$ complex number corresponds to $|s|$-order derivative operator of $f$. The results are applied for f being a monom. And remarkable relations concerning the e and $\pi$ order integrals and e and $\pi$ order derivatives are given, for e and $\pi$ transcendental numbers. Similar results may also be obtained for anticausal functions. For particular values of $a$ and $s$, the operator gives exactly Liouville fractional integral, Riemann fractional integral, Caputo fractional derivative, Liouville-Caputo fractional derivative. Finally, the new operator is neither integral nor derivative operator. It is integral and derivative operators at the same time. It deserves of being named : raoelinian operator is proposed.*

## I. INTRODUCTION

We tackle the problem of integral and derivative using the operator approach. And we are looking for the unique operator giving integrals and derivatives at the same time and for any order.

Our guidelines of thought are the following.

**I.1** Is it possible to derive the integral operator from the derivative operator or inversely ?

We have shown that it is easier to derive derivative operator from integral operator rather than the inverse [1] [2]

**I.2** In classical definition, the order $s$ of the integral and derivative operators must be **positive integer numbers**. Is it possible to extend them to any real and complex numbers? We have shown that the answer is yes [1] [2]

**I.3** Our research has been based on the important remarks concerning trigonometric and exponential functions.

### I.3.1 *Trigonometric functions*

For positive integer $n$, it is easily shown by iterating $n$ times the derivative operators $D^1$ and integral operator $J^1$ according to the following results

$$\begin{cases} D^1(sin)(x) = sin\left(x + \frac{\pi}{2}\right) & D^n(sin)(x) = sin\left(x + n\frac{\pi}{2}\right) \\ D^1(cos)(x) = cos\left(x + \frac{\pi}{2}\right) & D^n(cos)(x) = cos\left(x + n\frac{\pi}{2}\right) \end{cases}$$

$$\begin{cases} J^1(sin)(x) = sin\left(x - \frac{\pi}{2}\right) & J^n(sin)(x) = sin\left(x - n\frac{\pi}{2}\right) \\ J^1(cos)(x) = cos\left(x - \frac{\pi}{2}\right) & J^n(cos)(x) = cos\left(x - n\frac{\pi}{2}\right) \end{cases}$$

In geometric representation, $D^1$ and $D^n$ correspond to the rotations of $\frac{\pi}{2}$ and $n\frac{\pi}{2}$ respectively. $J^1$ and $J^n$ correspond to the rotations of $-\frac{\pi}{2}$ and $-n\frac{\pi}{2}$ respectively.

Then, we may define one new operator $R$ such

$$R^n(sin)(x) = sin\left(x - n\frac{\pi}{2}\right)$$

for any $n \in \mathbb{R}$. For $n$ positive (respectively negative) integer, we obtain the ordinary integral (respect. derivative) operator. The extension for $n \in \mathbb{C}$ is also possible.

These relations are utilized in the Fresnel representation to study RLC circuits (passive resistor R, self-inductance L and capacitance C) for sinusoidal electricity in terms of phase difference. The potential difference $v$ for R is in phase with the current $i$ ($v = Ri$). $v$ is delayed in phase of $\frac{\pi}{2}$ with $i$ for C ($v = integral\ of\ \frac{i}{C}$). $v$ is faster of $\frac{\pi}{2}$ with $i$ for L ($v = L\frac{di}{dt}$).

### I.3.2 *Exponential function*

The second remark is related to exponential functions

$$\begin{cases} D^1(e^{kx}) = ke^{kx} \\ J^1(e^{kx}) = k^{-1}e^{kx} \end{cases} \qquad \begin{cases} D^n(e^{kx}) = k^n e^{kx} \\ J^1(e^{kx}) = k^{-n} e^{kx} \end{cases}$$

We may then define one new operator $R$ such as

$$R^n(e^{kx}) = k^{-n}(e^{kx})$$

$R$ corresponds to the integral (respect.derivative) operator for positive (respect. negative) values of $n$.

In algebraic representation, the derivative operator is expressed in terms of multiplication, the integral operator by a division. This remark is utilized to solve an integro–differential equation with constant coefficients by postulating an exponential solution $exp(rx)$ with $r$ constant. The integro-differential equation is transformed into an algebraic equation with constant coefficients. Then we have to look for roots $r$, a well known problem.

The result is also utilized to study RLC circuits by the so-called complex number method.

I.3.3 **Remark**: the two approaches I.3.1 and I.3.2 are equivalent by using the Euler's relation between exponential and trigonometric functions.

## II. DEFINITION OF THE OPERATOR $R^s$

### *Theorem*

Let $E$ be the set of integrable and derivable causal functions defined on the interval $I = [a, +\infty[$, $a \in \mathbb{R}$ such $f(x) = 0$ for $x \le a$. Then the operator $R^s$

$$R^s(f)(x) = \frac{1}{\Gamma(s)} \int_a^x (x-y)^{s-1} f(y) dy$$

$$= \frac{x^s}{\Gamma(s)} \int_{\frac{a}{x}}^1 (1-u)^{s-1} f(ux) du$$

where $\Gamma(s)$ is the extension of Euler's gamma function for any $s$ ($s \in \mathbb{R}, s \in \mathbb{C}$) gives at the same time the integrals and derivatives operators at any order $s$. If $s \in \mathbb{R}_+$or $\mathcal{R}e(s) > 0$ ($s \in \mathbb{C}$), the operator $R^s$ gives the extension of the integral operator $J^s$of $f$ at any order $s$. If $s \in \mathbb{R}_-$or $\mathcal{R}e(s) < 0$, the operator $R^s$ corresponds to the extension of derivative operator $D^s$ at any order $|s|$.

### *Proof*

The proof has been given in our paper [1] for real order $s$, in our paper [2] for complex order $s$, in our paper [3] for real and complex order for causal functions and for real and complex order for anticausal functions [4]. Then $R^s$ contains the integral and derivative operators at any order $s$.
Let us recall the unique strategy applied in our approach [1] [2][3] [4]

1) Define the integral operator $J^1(f)(x)$ and the derivative $D^1(f)(x)$ of a causal functions $f$ at first order $s$.
2) Define $J^s(f)(x)$ and $D^s(f)(x)$ for any $s \in \mathbb{N}$.
3) Extend to $s \in \mathbb{Z}$, then to $s \in \mathbb{R}$
4) Extend to $s \in \mathbb{C}$
5) and finally, look for in which case we have the following conditions
  5a) the principle of correspondence

$$\lim_{s \to n} J^s(f)(x) = \int_a^x \int_a^{t_1} \int_a^{t_2} \dots \int_a^{t_n} f(t_n)\, dt_n dt_{n-1} \dots dt_1$$

for $n \in \mathbb{N}$

$$\lim_{s \to n} D^s(f)(x) = \frac{d^n}{dx^n}(f)(x)$$

for $n \in \mathbb{N}$. $\frac{d^n}{dx^n}$ is the ordinary derivative operator for $n$-order

5b) linear property of $J^s$ and $D^s$.

## III. PROPERTIES OF THE OPERATOR $R^s$

### III.1 **Linear property**

It is easy to demonstrate the relation

$$R^s(\alpha f + \beta g)(x) = \alpha R^s(f)(x) + \beta R^s(g)(x)$$

for any $f \in E, g \in E$. $\alpha \in \mathbb{C}$, $\beta \in \mathbb{C}$.

### III.2 **Semi-group property of $R^s$**

$$R^{s_1} R^{s_2} = R^{s_1+s_2} = R^{s_2} R^{s_1}$$

for any $s_1$ and $s_2$. It is assumed that the derivative operator $D^1$ (respectively $J^1$) is the inverse of the integral operator $J^1$ (respectively $D^1$).

### *Proof*

The proof is given in Appendix .
This property is very useful because it simplifies the demonstration of many relations in particular to define the derivative operator (see section IV.2)

### III.3 **Principle of correspondence**

It is easily shown the following properties

$$\lim_{s \to n} R^s(f)(x) = R^n(f)(x) = J^n(f)(x) \quad \forall n \in \mathbb{N}, \forall f \in E$$

$$\lim_{s \to -m} R^s(f)(x) = R^{-m}(f)(x) = D^m(f)(x) \quad \forall m \in \mathbb{N}, \forall f \in E$$

where $J^n$ is the ordinary integral of order $n$, and $D^m$ is the ordinary derivative operator of order $m$.

## IV. OBTENTION OF THE INTEGRAL OPERATOR $J^n$ AND DERIVATIVE OPERATOR $D^m$ FROM $R^s$

### IV.1 **Obtention of the integral operator $J^n$ from $R^s$**

It may be easily shown

$$R^n(f)(x) = J^n(f)(x) \quad \forall n \in \mathbb{R}_+ \text{ or } \mathcal{R}e(s) > 0 \text{ if } n \in \mathbb{C}$$

Then

$$R^n = J^n$$

### IV.2 **Obtention of the derivative operator from $R^s$**

We have two possibilities to define derivative operators from $R^s$.

#### IV.2.1 *The left-hand derivative operator $D_L^m$*

$$D_L^m = R^{-k} R^{k-m} = D^k R^{k-m}$$

for any $k \in \mathbb{N}$ and any $m \in \mathbb{R}_+$ or $\mathcal{R}e(s) > 0$ if $m \in \mathbb{C}$ no sommation on $k$

or

$$D_L^m(f)(x) = D^k R^{k-m}(f)(x)$$

$$= D^k \frac{1}{\Gamma(k-m)} \int_a^x (x-y)^{k-m-1} f(y) dy$$

By successive application of derivative operator $D$ on the integral, it is easily found that $D_L^m(f)(x)$ is independent on $k$ and we obtain

$$D_L^m(f)(x) = D^1 \frac{1}{\Gamma(1-m)} \int_a^x (x-y)^{-m} f(y) dy$$

If $f(x) = C$ where $C$ is a constant

$$D_L^m(C) = C \frac{(x-a)^{-m}}{\Gamma(1-m)} \neq 0 \ \text{if } C \neq 0$$

The Riemann fractional derivative is defined by [5]

$$D_{Ri}^m(f)(x) = \frac{d}{dx} \frac{1}{\Gamma(1-m)} \int_0^x (x-y)^{-m} f(y) dy$$

which is exactly $D_L^m(f)(x)$ for $a = 0$.

### IV.2.2 The right-hand derivative operator $D_R^m$

$$D_R^m = R^{k-m} R^{-k} = R^{k-m} D^k$$

for any $k \in \mathbb{N}$ and any $m \in \mathbb{R}_+$ or $\mathcal{R}e(s) > 0$ if $m \in \mathbb{C}$ no sommation on $k$.

$$D_R^m(f)(x) = R^{k-m}(D^k f)(x)$$
$$= \frac{1}{\Gamma(k-m)} \int_a^x (x-y)^{k-m-1} f^{(k)}(y) dy$$

By successive integration by part of the function under the integration sign, we obtain that the second member is independent on $k$.

$$D_R^m(f)(x) = \frac{1}{\Gamma(1-m)} \int_a^x (x-y)^{-m} f^{(1)}(y) dy$$

If $f(x) = C$ where $C$ is a constant

$$D_R^m(C) = 0$$

even if the constant $C \neq 0$

$D_R^m(f)$ is then the good choice instead of $D_L^m$ if you require the principle of correspondence that the derivative of a constant is null.

Let us prove now that

$$D_R^m(f)(x) = R^{-m}(f)(x)$$

We have

$$R^s(f)(x) = \frac{1}{\Gamma(s)} \int_a^x (x-y)^{s-1} f(y) dy$$

Let us integrate by part the function under the integration sign.

The integrated terms are null and we have

$$R^s(f)(x) = \frac{1}{\Gamma(s+1)} \int_a^x (x-y)^s f^{(1)}(y) dy$$

Let us change $s$ to $-m$

$$R^{-m}(f)(x) = \frac{1}{\Gamma(1-m)} \int_a^x (x-y)^{-m} f^{(1)}(y) dy$$
$$= D_R^m(f)(x)$$

then

$$R^{-m} = D_R^m$$

The Liouville-Caputo fractional derivative definition is [5]

$$_{LC}D_+^m(f)(x) = \frac{1}{\Gamma(1-m)} \int_{-\infty}^x (x-y)^{-m} f^{(1)}(y) dy$$

which is exactly $D_R^m(f)(x)$ with the particular value $a = -\infty$. The Caputo fractional derivative is defined by [5]

$$_C D_+^m(f)(x) = \frac{1}{\Gamma(1-m)} \int_0^x (x-y)^{-m} [\frac{d}{dy} f(y)] dy$$

which is $D_R^m(f)(x)$ with the particular value $a = 0$.

### IV.2.3 Remark

If the derivative operator $D^1$ (respectively $J^1$) is the inverse of the integral operator $J^1$ (respectively $D^1$) then we have the semi-group property for the operator $R^s$ and we have one derivative operator $D_L^m = D_R^m$. If it is not the case we don't have the semi-group property for the operator $R^s$ but the semi-group property stands for $J^s$ and $D^s$.

## V. APPLICATIONS OF THE RESULTS [1]

Let us apply now our results to the function $f(x) = x^p$ for any positive real $p$.

$$R^s(x^p) = \frac{\Gamma(p+1)}{\Gamma(p+s+1)} x^{p+s} = \frac{p!}{(p+s)!} x^{p+s}$$

for any positive real numbers $s$ and $p$.

$$J^\pi(x^e) = R^\pi(x^e) = \frac{\Gamma(e+1)}{\Gamma(e+\pi+1)} x^{\pi+e} = \frac{e!}{(e+\pi)!} x^{\pi+e}$$

$$J^e(x^\pi) = R^e(x^\pi) = \frac{\Gamma(\pi+1)}{\Gamma(\pi+e+1)} x^{e+\pi} = \frac{\pi!}{(e+\pi)!} x^{\pi+e}$$

The ratio

$$\frac{R^\pi(x^e)}{R^e(x^\pi)} = \frac{e!}{\pi!}$$
$$= 0.592\ 761\ 747\ 048\ 502\ 880\ 285\ 354\ 552\ 437\ 32$$

is independent on $x$.

We obtain easily the derivative operator $D_R^e$ from $R^{-e}$

$$D_R^e(x^\pi) = R^{-e}(x^\pi) = \frac{\Gamma(\pi+1)}{\Gamma(\pi-e+1)} x^{\pi-e} = \frac{\pi!}{(\pi-e)!} x^{\pi-e}$$

$$D_R^\pi(x^e) = R^\pi(x^e) = \frac{\Gamma(e+1)}{\Gamma(e-\pi+1)} x^{e-\pi} = \frac{e!}{(e-\pi)!} x^{e-\pi}$$

Then the product

$$D_R^e(x^\pi) . D_R^\pi(x^e) = \frac{e! \, \pi!}{(\pi - e)! \, (e - \pi)!}$$

$$= 22.364\ 994\ 517\ 058\ 857\ 454\ 906\ 921\ 720\ 114$$

is independent on $x$.
Direct calculation have been given in our paper [1].

## VI. CONCLUSIONS

In 1695 Leibnitz raised for the first time the problem of fractional derivatives: "Can the meaning of derivatives with integer order be generalized to derivative with noninteger orders?" L'Hospital replied to Leibnitz by another question "what if the order is ½" and Leibnitz answered "It will lead to a paradox from which one day useful consequences will be drawn".

Several approaches have been done and detailed bibliography may be found in references [7] and [8].

We think that the question of Leibnitz has now fully received a satisfactory answer not only for derivatives but for integrals too because the order is extended to any real numbers and complex numbers for derivatives and integral orders.

The problem may be considered as the definitions of the s-power of an operator for any real $s$ and any complex $s$ too but we will not develop this point of view.

The operator $R$ is a new one, it is neither integral operator nor derivative operator, it is the generalization of integral and derivatives operators at the same time. It has many interesting properties. We think that it deserves then of being named: raoelinian operator is proposed.

## APPENDIX

We have to show

$$R^{s_1} R^{s_2}(f)(x) = R^{s_1+s_2}(f)(x) = R^{s_2} R^{s_1}(f)(x)$$

for any $f \in E$ and for any $s_1$ and $s_2$. It is assumed that the operator $J^1$ is the inverse of the operator $D^1$. If it is not the case then we don't the semi-group property for $R^s$.

$$R^{s_1} R^{s_2}(f)(x) =$$
$$\frac{1}{\Gamma(s_1)\Gamma(s_2)} \int_a^x dy\, (x - y)^{s_1-1} \int_a^y dz\, (y - z)^{s_2-1} f(z)$$

By application of Dirichlet's formula given by Whittaker and Watson [5] [6]

$$\int_a^x dy\, (x - y)^{s_1-1} \int_a^y dz\, (y - z)^{s_2-1} g(y, z)$$

$$= \int_a^z dz \int_z^x dy\, (x - y)^{s_1-1}(y - z)^{s_2-1} g(y, z)$$

For $g(y, z) = f(z)$, we obtain

$$R^{s_1} R^{s_2}(f)(x)$$
$$= \frac{1}{\Gamma(s_1)\Gamma(s_2)} \int_a^x dz f(z) \int_x^z dy\, (x - y)^{s_1-1}(y - z)^{s_2-1}$$

we perform the change of variable

$$u = \frac{y - z}{x - z}$$

$$y = z + u(x - z) \qquad dy = (x - z) du$$

$$R^{s_1} R^{s_2}(f)(x)$$

$$= \frac{1}{\Gamma(s_1)\Gamma(s_2)} \int_a^x dz f(z)\, (x - z)^{s_1+s_2-1} \int_0^1 du\, (1 - u)^{s_1-1} u^{s_2-1}$$

we obtain the Euler's beta function

$$B(s_1)(s_2) = \frac{\Gamma(s_1)\Gamma(s_2)}{\Gamma(s_1 + s_2)}$$

$$R^{s_1} R^{s_2}(f)(x) = \frac{1}{\Gamma(s_1 + s_2)} \int_a^x dz f(z)\, (x - z)^{s_1+s_2-1}$$
$$= R^{s_1+s_2}(f)(x)$$

$R^{s_1} R^{s_2}$ is symmetric in $(s_1, s_2)$, then

$$R^{s_1} R^{s_2}(f)(x) = R^{s_1+s_2}(f)(x) = R^{s_2} R^{s_1}(f)(x)$$

We obtain the semi-group property of $R^s$